\documentclass[a4paper]{article}

\usepackage{latexsym}
\usepackage{graphicx}
\usepackage{amsmath}
\usepackage{amssymb}
\usepackage{amscd}
\usepackage{epsf}

\newcommand{\lbi}{\ensuremath{L_B(\infty)}}

\newcommand{\lb}{\ensuremath{L_B}}

\newcommand{\spmi}{\ensuremath{SPM(\infty)}}
\newcommand{\spm}{\textsc{SPM}}

\newcommand{\proof}{\vspace*{-0.2cm} \textbf{Proof : }}
\newcommand{\cqfd}{\hfill \fbox{} \vskip 0.2cm}

\newcommand{\flips}{\ensuremath{\stackrel{*}{\rightarrow}}}

\newcommand{\centre}[1]{\raisebox{-1ex}{#1}}

\newcommand{\ran}[2]{\ensuremath{#1 \rightarrow #2}}

\newcommand{\fig}[2]{\begin{figure}[!h] \centerline{\includegraphics{#1.eps}}
\caption{\label{fig_#1} #2} \end{figure}}

\begin{document}

\newtheorem{prop}{Proposition}
\newtheorem{theo}{Theorem}
\newtheorem{corol}{Corollary}
\newtheorem{lemme}{Lemma}
\newtheorem{definition}{Definition}
\newtheorem{notation}{Notation}
\newtheorem{rem}{Remark}

\begin{center}
{\huge \textbf{Generalized Integer Partitions,\\ Tilings of Zonotopes and Lattices}}\\
\vskip 0.2cm
{\large{Matthieu Latapy\,\footnote{\textsc{liafa}, Universit\'e Paris 7, 2 place Jussieu, 75005 Paris. latapy@liafa.jussieu.fr}}}
\vskip 0.2cm
\end{center}

\textbf{Abstract :} In this paper, we study two kinds of combinatorial
objects, generalized integer partitions and tilings of two dimensional
zonotopes, using
dynamical systems and order theory. We show that the sets of partitions
ordered with a simple dynamics, have the distributive lattice structure.
Likewise, we show that the set of tilings of zonotopes, ordered with
a simple and classical dynamics, is the disjoint union of distributive
lattices which we describe.
We also discuss the special case of linear integer
partitions, for which other dynamical systems exist.
These results give a better understanding of the behaviour of
tilings of zonotopes with flips and dynamical systems
involving partitions.

\vskip 0.2cm
\textbf{Keywords :} Integer partitions, Tilings of Zonotopes,
Random tilings, Lattices, Sand Pile Model, Discrete Dynamical Systems.
\vskip 0.5cm

\section{Preliminaries}

In this paper, we mainly deal with two kinds of combinatorial objects: integer
partitions and tilings. An \emph{integer partition problem} is a set of
partially ordered variables, and a \emph{partition} is an affectation of an
integer to each of these variables such that the order between these values
is compatible with the order over the variables.
Tilings are coverings of a given space with
tiles from a fixed set. We are concerned here with tilings of
two dimensional zonotopes (i.e. $2D$-gons) with lozenges, usually
called \emph{random tilings} in theoretical physics.
These two apparently different kinds of objects have been brought
together in \cite{Des97,DMB99}. Some special cases were already known
from \cite{DMB97, Bai97, Bai99} but the correlations between partitions
and tilings are treated in general for the first time in \cite{Des97}.
We begin with a study of the structure of the set of solutions
to a partition problem, and then use the obtained results and the
correspondence between partitions and tilings to
study tilings of zonotopes with flips.

We mainly use two tools: dynamical systems and orders. The use of
\emph{dynamical systems} allows an intuitive presentation of the
results and makes it easier to understand the relations between
the different concerned objects. The use of orders is natural since
they appear as structures of the sets we study. An \emph{order relation}
is a binary relation over a set, such that for all $x$, $y$ and $z$ in
this set, $x$ is in relation with itself (reflexivity), the fact that
$x$ is in relation with $y$ and $y$ is in relation with $z$ implies
that $x$ is in relation with $z$ (transitivity),
and the fact that $x$ is in relation with $y$ and $y$ with $x$
implies $x=y$ (antisymmetry). The set is then a \emph{partially
ordered set} or, for short, a \emph{poset}.
Now, in a poset, if any two elements have an \emph{infimum},
i.e. a greatest lower element, and a \emph{supremum}, i.e. a lowest
greater element, then it is a \emph{lattice}.
The infimum of two elements $a$ and $b$ in a lattice $L$ is denoted
by $\inf_L(a,b)$,
and their supremum is denoted by $\sup_L(a,b)$. We often
write simply $\inf(a,b)$ and $\sup(a,b)$ when the context makes
it clear which lattice is concerned.
A lattice is \emph{distributive} if for all $a$, $b$ and $c$:
$\sup(\inf(a,b),\inf(a,c)) = \inf(a,\sup(b,c))$ and 
$\inf(\sup(a,b),\sup(a,c)) = \sup(a,\inf(b,c))$.
For more details, see
\cite{DP90}.

\section{The lattices of integers partitions}
\label{sec_part}

\subsection{Generalized integer partitions}
\label{sec_gl}

A \emph{generalized integer partition problem}, or simply a 
\emph{partition problem} \cite{Sta72, DMB97}, is defined by a Directed Acyclic
Graph (DAG) $G = (V,E)$ and a positive integer $h$. A solution
of such a problem, called a \emph{partition}, is a function
$a : V \rightarrow \lbrace 0,1,\dots,h\rbrace$ such that
$a(v) \ge a(w)$ for all $v$ and $w$ in $V$
      such that there is an edge from $v$ to $w$ in $G$.
The integer $a(v)$ is usually denoted by $a_v$, and the set of all
the solutions of a given partition problem $(G,h)$ is $P(G,h)$.
The graph $G$ is called the \emph{base} of the partition problem,
and $h$ is the height of the problem. We extend here the usual definition
by allowing $h = \infty$, which means that the parts are unbounded.
For example, if $G_1=(V_1,E_1)$ is defined by
$V_1 = \lbrace v_1,v_2,v_3,v_4,v_5 \rbrace$
and $E_1 = \lbrace (v_1,v_2) , (v_1,v_3), (v_2,v_4), (v_3,v_4), (v_4,v_5) \rbrace$ then a
possible solution to the partition problem $(G_1,3)$ is
$(v_1\mapsto 3, v_2\mapsto 3, v_3\mapsto 2, v_4\mapsto 1, v_5\mapsto 0)$.


We can obtain all the solutions of a partition problem
$(G=(V,E),h)$ with the following dynamical system.
We consider a partition $a$ as a \emph{state} of the system, and $a_v$
as a number of grains stored at $v \in V$. The initial state of the
system is the empty one: $a_v=0$ for all $v\in V$.
Then, the following rule is iterated:\\
\centerline{\emph{A grain can be added at $v \in V$ if the
obtained configuration is a partition.}}\\
In other words, if $a$ is the state of the system, then the transition
$a \longrightarrow b$ is possible if $\forall \nu\not= v$,
$b_\nu=a_\nu$, $b_v=a_v+1$ and $b$ is a partition.
The state $b$ is then called a \emph{successor} of $a$.
See for example Figure~\ref{fig_ex_dyn}.

\fig{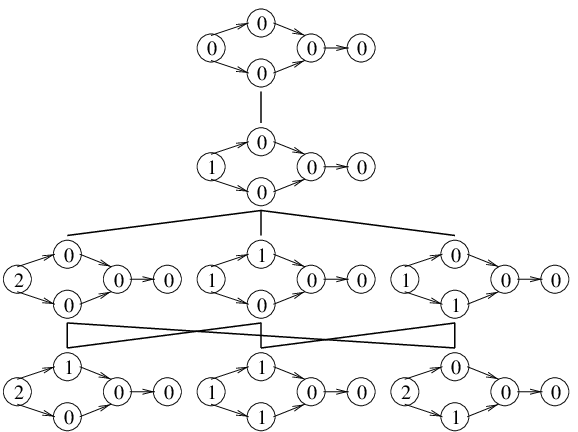}{Example of the dynamics with the partition problem $(G_1,2)$:
the first partitions of the set $P(G_1,2)$ with the possible transitions.
We show in each vertex the number of grains it contains. Notice that we do
not display all the reachable partitions in this diagram, but only the
first ones.}

Since there can be no cycle in a sequence of states of the system,
this transition rule induces a partial order over the possible states, i.e.
over the set of partitions $P(G,h)$. We can now state the main
result of this section, which will be useful during the study
of tilings problems in Section~\ref{sec_til}.

\begin{theo}
\label{th_part}
Given a partition problem $(G=(V,E),h)$, the set $P(G,h)$
equipped with the order induced by the transition rule is a
\emph{distributive lattice}.
Moreover, the infimum (resp. supremum) of two given partitions
$a$ and $b$ in this set is the partition $c$ (resp. $d$) defined by:
$$\forall v \in V,\  c_v=\mbox{max}(a_v,b_v),$$
$$\forall v \in V,\  d_v=\mbox{min}(a_v,b_v).$$
\end{theo}
\proof
It is clear that $c$ is a partition. Consider now a partition $\gamma$. If
for one vertex $\nu \in V$, $\gamma_{\nu} < c_{\nu}=\max(a_{\nu},b_{\nu})$
then $\gamma$ is clearly unreachable from $a$ or $b$ by iteration of the
transition rule. Likewise, if for all $\nu \in V$,
$\gamma_{\nu} > c_{\nu}$ then $\gamma$ is
clearly reachable from $c$. Therefore, we have that $c=\inf(a,b)$.
The proof for $d=\sup(a,b)$ is similar. Therefore, $P(G,h)$ is a lattice.
Now, it is easy from these formula to verify that the properties
required for a lattice to be distributive are fulfilled.
\cqfd

This result says for example that the set $P(G_1,2)$, partially
shown in Figure~\ref{fig_ex_dyn}, is a distributive lattice.
If $h=\infty$,
we obtain an infinite lattice which contains all the possible
partitions over the base graph of the problem. Moreover, it is
easy to verify that the sets $P(G,h)$ with $h < \infty$ are
sub-lattices of the infinite one.

Let us recall that an order \emph{ideal} $I$ is a subset
of an order $P$ such that $x \in I$ and $y \le x$ in $P$ implies
$y \in I$. In the finite case, an ideal is defined by a set
of uncomparable elements of $P$ and contains all the elements of $P$
lower than any element of this set \cite{DP90}.
From \cite{DP90}, we know that for any distributive lattice
there exists a unique order such that the given lattice is isomorphic
to the lattice of the ideals of the order, ordered by inclusion.
Conversely, the set of ideals such ordered is always a distributive
lattice. We will now show that for any partition problem $(G,h)$
with $h < \infty$, $P(G,h)$ is isomorphic to the lattice of
ideals of a partial order, which is another way to prove and
understand Theorem~\ref{th_part}.

First notice that, since $G=(V,E)$ is a Directed Acyclic Graph,
it can be viewed as an
order over the vertices of $G$. Now consider the ordered set
$\lbrace 1,2,\dots,h \rbrace$ with the natural order
$1 < 2 < \dots < h$ and the direct product 
$G \times \lbrace 1,2,\dots,h \rbrace$ defined by:
for all $a$, $b$ in $V$, $i$, $j$ in $\lbrace 1,2,\dots,h \rbrace$,
$(a,i) \le (b,j)$ if and only if $a \le b$ in the order induced
by  $G$ and $i \le j$. See Figure~\ref{fig_ideals} for an example.
Consider now an ideal $I$ of $G \times \lbrace 1,2,\dots,h \rbrace$.
We can define the application $p_I$, from the set $V$ of vertices of $G$
to $\lbrace 0,1,2,\dots,h \rbrace$ by $p_I(v) = j$ where $j$ is the
maximal integer $i$ such that $(v,i) \in I$, if any, and
$p_I(v) = 0$ if there exists no such $i$.
It is clear that the application defined this way is a partition,
and that $r : I \mapsto p_I$ is a bijection between the set of
ideals of $G \times \lbrace 1,2,\dots,h \rbrace$ and the set
of partitions $P(G,h)$.
Now, we will see that $r$ is actually an order isomorphism by showing
that it preserves the covering relation (i.e. the transitive reduction
of the order relation) of these orders. Recall that a partition $p$
in $P(G,h)$ is covered exactly by the partitions $p'\in P(G,H)$ such that
$p'$ has one more grain at one vertex, say $v$. This implies that all the
vertices $v'$ such that there is a path from $v'$ to $v$ in $G$ must have
strictly more grains than $v$. This means that the corresponding ideals
$I = r^{-1}(p)$ and $I' = r^{-1}(p')$ verify
$I' \setminus I = \lbrace (v,p_v+1) \rbrace$, which is exactly the
covering relation in the lattice of ideals of
$G \times \lbrace 1,2,\dots,h \rbrace$. Conversely, let us consider
two orders ideals $I$ and $I'$ such that
$I' \setminus I = \lbrace (v,i) \rbrace$. Then, it is clear that
the corresponding partitions $r(I)$ and $r(I')$ only differ by $1$
at vertex $v$.



\fig{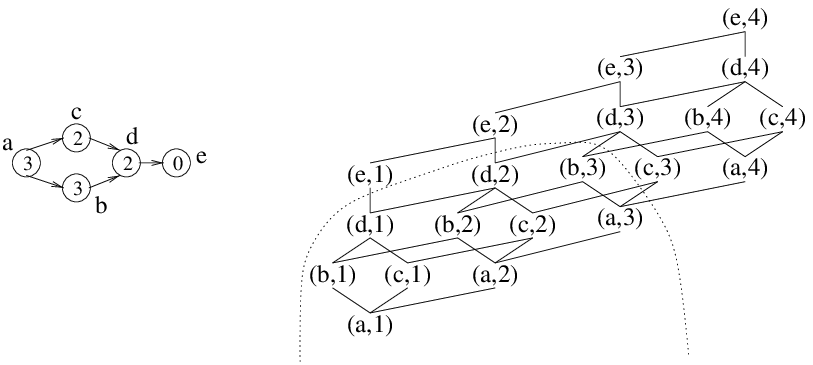}{The solutions of $(G,h)$ are nothing but the ideals
of $G \times \lbrace 1,2,\dots,h \rbrace$. Here, we show a graph $G$
(left) with vertices $\lbrace a,b,c,d,e\rbrace$ and the product
$G \times \lbrace 1,2,3,4\rbrace$ (right), which is equivalent
to the partition
problem $(G,4)$. We show the ideal equivalent to the
partition displayed on the graph.}

We will now use Theorem~\ref{th_part} to
study special classes of partition problems, the hypersolid ones,
and have a special attention for the so-called linear partitions.
However, the reader mostly interested in tilings may directly go to
Section~\ref{sec_til}.

\subsection{Hypersolid and Linear partitions}

When the base graph of a partition problem is a directed
d-dimensional grid then the problem is called \emph{hypersolid}
and the solutions are called \emph{hypersolid partitions} \cite{Mac16}.
Formally, a $d$-dimensional grid is a graph $G=(V,E)$
with $V = (\mathbb{N}^*)^d$ such that there is a path from $v=(v_1,v_2,\dots,v_d)$ to
$w=(w_1,w_2,\dots,w_d)$ if and only if $\forall i, v_i \le w_i$.
For example, the $1$-dimensional hypersolid partition problem
of length $l$ and height $h$ is $(G=(V,E),h)$ where $V=\lbrace
1,2,\dots ,l \rbrace$ and $E=\lbrace(i,i+1)|1\le i<l\rbrace$.
A $1$-dimensional hypersolid partition is called a \emph{linear
partition}, and it is denoted by $a=(a_1,a_2,\dots,a_l)$.
A $2$-dimensional hypersolid partition is called a \emph{plane
partition} and is usually described by an array such that the position
$(i,j)$ contains $a_{i,j}$. For example, 
\scalebox{.5}{\begin{tabular}{|c|c|c|} \hline 4&4&2\\ \hline 3&2&2\\ \hline 1&1&0\\ \hline \end{tabular}}
is a solution of the plane partition problem of size $3\times 3$ and of
height $5$.
In this section, we will denote by $H(d,s,h)$ the $d$-dimensional
hypersolid partition problem of size $s$ and height $h$.
Likewise, we denote by $P(H(d,s,h))$ the set of the solutions
of $H(d,s,h)$, ordered by the reflexive and transitive closure of the
transition relation described in Section~\ref{sec_gl}.
Therefore, $P(H(1,\infty,\infty))$ is the set of all
the linear partitions.

The hypersolid partitions are of particular interest and have been widely
studied \cite{And76}. Theorem~\ref{th_part} tells us that the sets
$P(H(d,s,h))$ are distributive lattices.
This is a result which is also known as a
consequence of the bijection between $P(H(d,\infty,\infty)$ and
the lattice of ideals of $\mathbb{N}^d$ ordered by inclusion
\cite{DP90}.


We will now consider the special case of linear partitions. They have
been widely studied as a fundamental combinatorial object \cite{And76}.
As mentioned above, a linear partition of an integer $n$ is simply
a decreasing sequence of integers,
called \emph{parts}, such that the sum of the parts is exactly $n$.
A linear partition is usually represented by its Ferrer's diagram,
a sequence of columns
such that if the $i$-th part is equal to $k$ then the $i$-th
column contains exactly $k$ stacked squares, called grains.
In 1973, Brylawski proposed a dynamical system to study these partitions
\cite{Bry73}:
given a partition $a$, a grain can fall from column $i$ to
column $i+1$ if $a_i - a_{i+1} \ge 2$ and a grain can slip from
column $i$ to column $j>i+1$ if for all $i<k<j$, $a_k = a_i-1 = a_j+1$. See
Figure~\ref{fig_regles}.

\fig{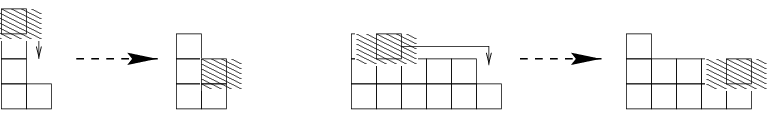}{The two evolution rules of the dynamical system defined by Brylawski.}

Brylawski showed that the iteration of these rules
from the partition $(n)$ gives the lattice of all the linear partitions
of $n$ ordered with respect to the dominance ordering defined by:
$$a \ge b \mbox{ if and only if }
\sum_{i=1}^j a_i \ge \sum_{i=1}^j b_i \mbox{ for all $j$,}
$$
i.e. the prefix sums of $a$ are greater than or equal to
the prefix sums of $b$. This lattice is denoted by $L_B(n)$. See
Figure~\ref{fig_diag} (left) for an example.
If one iterates only the first rule defined by Brylawski, one obtains
the Sand Pile Model and the set of linear partitions obtained
from $(n)$ is a lattice, denoted by $SPM(n)$, with respect to the dominance
ordering  \cite{GK93}. See Figure~\ref{fig_diag} (right) for an example.
In \cite{LMMP98} and \cite{LP99}, it is proved that when these systems
are started with one infinite first column
the sets of reachable configurations have a structure of
infinite lattice, denoted by $SPM(\infty)$ and $L_B(\infty)$.
It is also shown in these papers that, if we consider $a$ and $b$ in
$SPM(\infty)$ or $L_B(\infty)$ then their infimum $c$ is defined by:
\begin{align}
c_i = max(\sum_{j\geq i} a_j, \sum_{j \geq i} b_j) - \sum_{j >i} c_j\hskip
0.3cm \mbox{ for all } i
\label{formule_inf}
\end{align}
The lattice $L_B(\infty)$ contains all
the linear partitions, just as $P(H(1,\infty,\infty))$.
We will now study the connection between the dynamical system
defined by Brylawski and the one defined in Section~\ref{sec_gl}.

\fig{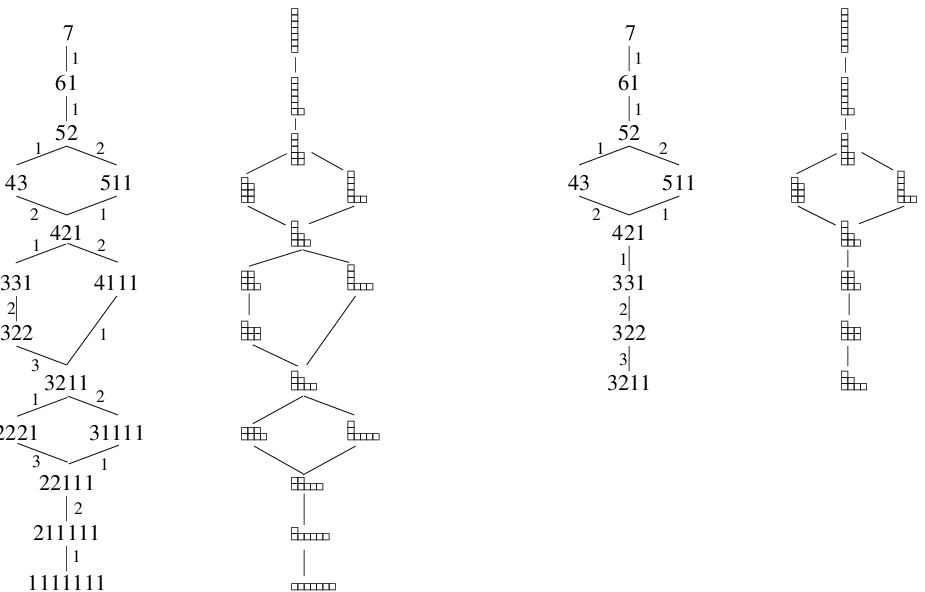}{On the left, the diagram of the lattice $L_B(7)$, and on the
right the diagram of $SPM(7)$. We showed the representation by
piles of grains, and we displayed on each edge the column
from which the grain falls during the corresponding transition.}

\begin{theo}
The application:
$$
\begin{array}{lclc}
\pi_{\lb} : & \lbi & \longrightarrow & P(H(1,\infty,\infty))
\end{array}
$$
such that $\pi_{\lb}(a)_i$ is equal to $\sum_{j\ge i}a_j$
is an order embedding which preserves the infimum.
\end{theo}
\proof
To clarify the notations, let us denote by $\pi$ the application
$\pi_{\lb}$ in this proof.
Let $a$ and $b$ be two elements of \lbi. We must show that $\pi(a)$
and $\pi(b)$ belong to $P(H(1,\infty,\infty))$, that $a \ge_{\lbi} b$ is
equivalent to $\pi(a) \ge_{P(H(1,\infty,\infty))} \pi(b)$ and that
$\inf_{P(H(1,\infty,\infty))}(\pi(a),\pi(b)) = \pi(\inf_{\lbi}(a,b))$.
The two first points are easy:
$\pi(x)$ is obviously a decreasing sequence of integers for any $x$, and
the order is preserved.
Now, let $u = \inf(a,b)$. Then,
$$
\begin{array}{lclr}
\pi(u)_i & = & \sum_{j\ge i}u_j &\\
         & = & max(\sum_{j\ge i}a_j,\sum_{j\ge i}b_j)& \mbox{from~(\ref{formule_inf})}\\
         & = & max(\pi(a)_i,\pi(b)_i) \\
         & = & \inf(\pi(a),\pi(b))_i
\end{array}
$$
which proves the claim.
\cqfd

Notice that if we consider the restriction of $\pi_{\lb}$ to
\spmi, denoted by $\pi_{\spm}$, a similar proof shows that
$\pi_{\spm}$ is an order embedding which preserves the infimum.
However, these orders embeddings are not \emph{lattices} embeddings, since
they do not preserve the supremum. For example, if $a=(2,2)$ and $b=(1,1,1)$,
then $\pi_{\lb}(a) = (4,2)$, $\pi_{\lb}(b) = (3,2,1)$,
$c = \sup_{\lbi}(a,b) = (2,1)$ but $\pi_{\lb}(c) = (3,1)$ and 
$\sup_{P(H(1,\infty,\infty))}((4,2),(3,2,1)) = (3,2)$.
Notice that there can be no lattice embedding from
$\lbi$ to $P(H(1,\infty,\infty)$ since the fact that $P(H(1,\infty,\infty)$
is a \emph{distributive} lattice would imply that $\lbi$ would
be distributive, which is not true.


\section{The lattice of the tilings of a zonotope}
\label{sec_til}

A \emph{tiling problem} (see Lecture 7 in \cite{Zie95} for example)
is defined by a finite
set of tiles $T$, called the prototiles, and a polygon $P$.
A solution of the problem is
a \emph{tiling}: an arrangement of translated copies of prototiles
which covers exactly $P$ with no gap and no overlap.
The study of these problems is a classical field in mathematics.
They appear in computer science with the famous result of Berger
\cite{Ber66}, who proved the undecidability of the problem of
knowing if, given a set of prototiles, the whole plane can be
tiled using only copies of prototiles.


We are concerned here with tilings of $d$-dimensional zonotopes
with rhombic tiles. A \emph{$d$-dimensional zonotope} $Z$ is
defined from a family of $d$-dimensional vectors $\lbrace v_1,
v_2, \dots, v_D \rbrace$ and a set of positive integers $\lbrace l_1,
l_2, \dots, l_D \rbrace$ by:
$$
Z = \lbrace \sum_{i=1}^D \alpha_i v_i, 0 \le \alpha_i \le l_i \rbrace.
$$
If $l_i=1$ for all $i$, it is also called the \emph{Minkowski} sum
of the vectors. Notice that a two dimensional zonotope is a $2D$-gon
(lozenge, hexagon, octogon, etc).
The \emph{rhombic tiles} are obtained as the Minkowski sum of $d$ vectors
among the ones which generate the zonotope we want to tile. Such a tile
is called a \emph{rhombus}. If $d=2$, they are simply lozenges.
Rhombic tilings of zonotopes can be seen as projections of
sets of faces of a $D$-dimensional grid onto the $d$-dimensional
subspace along a generic direction.
By construction, the so-obtained tiles are the projections of the
$d$-dimensional facets of the $D$-dimensional grid and the tiled region
is a zonotope.
The tiling is then called a $D \rightarrow d$ \emph{tiling}
and the integer $D-d$ is the \emph{codimension} of the tiling.
Figure~\ref{fig_ex_til} shows an example of a \ran{3}{2}\ tiling
(left) and a \ran{4}{2}\ one (right). One can refer to
\cite{DMB97, DMB99, Des97, Bai97, Bai99} for more details and examples.

\fig{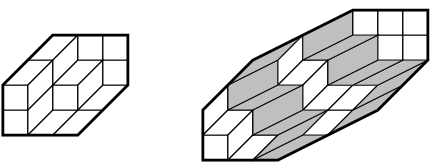}{Examples of tilings. Left : a $3\rightarrow 2$
tiling. Right : a $4\rightarrow 2$ tiling. Notice that if we delete
the shaded tiles in the $4\rightarrow 2$ tiling then we obtain
a $3\rightarrow 2$ tiling.}

A dynamical transformation is usually defined over \ran{D}{d}\ tilings
of a zonotope $Z$.
If we consider the zonotopes obtained as the Minkowski sum of $d+1$ vectors
among the ones which generate $Z$, then we obtain the most simple
sub-zonotopes\,\footnote{A sub-zonotope $Z'$ of a zonotope $Z$ generated
by $V=\lbrace v_1,\dots,v_D\rbrace$ and $\lbrace l_1,\dots,\l_D\rbrace$
is a zonotope generated by a subset $V'=\lbrace v_{i_1},\dots,v_{i_k}\rbrace$
of $V$ and $\lbrace l'_{i_1},\dots,\l'_{i_k}\rbrace$ with
$l'_{i_j} \le l_{i_j}$ for all $j$.}
of $Z$ with non-trivial tilings. Let us call
these zonotopes, generated by a family of $d+1$ vectors,
\emph{elementary zonotopes}. One can notice that there are exactly two ways
to tile such a zonotope, with $d+1$ tiles. For example, if $d=2$ then the
elementary zonotopes are hexagons and the possible tilings are
\centre{\includegraphics{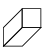}} and \centre{\includegraphics{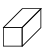}}. Therefore, we
can define a dynamics over the tilings of a given zonotope $Z$:
$t \longrightarrow t'$ if and only if we can obtain $t'$ by
changing in $t$ the tiling of one elementary sub-zonotope of $Z$.
This means that we locally rearrange $d+1$ tiles in the tiling $t$
to obtain $t'$.
In the $d=2$ case, it corresponds to the local rearrangement:
$\centre{\includegraphics{hexa1.eps}} \longrightarrow \centre{\includegraphics{hexa2.eps}}$.
If we call $t_1$ and $t_2$ the two possible tilings of an elementary
zonotope, we will call \emph{flip} the local transformation of $t_1$
into $t_2$ in any tiling $t$, and \emph{inverse} flip
the local transformation of $t_2$ into $t_1$ in any tiling $t$.
This gives an orientation to the notion of flips, and in the following
we will only be concerned with flips (not inverse ones),
unless explicitely specified. See Figure~\ref{fig_pav_tr} for some examples.

In order to continue with the relationship between tilings and
partitions, we need to recall the classical notions of
de Bruijn surfaces and families.
\emph{De Bruijn} grids \cite{Bru81,Bru86} are dual
representations of tilings which have been widely used to
obtain important results. De Bruijn grids are composed of
the de Bruijn $(d-1)$-dimensional surfaces, which join
together the middles of the two opposite sides of each
tile. Since the tiles are rhombus, it is always
possible to extend these surfaces through the tiling up to the
boundary. The set of tiles crossed by a de Bruijn surface is called a
\emph{worm}: it is composed of adjacent tiles. If $d=2$, de Bruijn
surfaces are lines, and they join together the middles of the opposite
edges of the lozenges tiles. An example is given in Figure~\ref{fig_bruijn}.
Each tile is crossed by exactly $d$ de Bruijn surfaces, and there
is no intersection of $d+1$ surfaces. On the other hand, there are
surfacs which can never intersect, even in an infinite tiling.
They join rhombus faces of same orientation, such as lines $a$
and $b$ in Figure~\ref{fig_bruijn}. We say that these surfaces belong
to the same family. A family is equivalent to an edge orientation.
In the following, we call \emph{de Bruijn family} the set of
tiles in the worms that correspond to the lines of a family.
To sum up, we can say that a de Bruijn family of tiles is defined
by a vector among the ones which generate the tiled zonotope,
and the family contains all the tiles which have this vector as an edge.
For example, in Figure~\ref{fig_ex_til} (right),
a de Bruijn family of tiles is shaded. Notice that
deleting such a family in a \ran{D}{d}\ tiling gives
a \ran{D-1}{d} tiling \cite{Eln97}. In the following, when we will
consider a \ran{D}{d}\ tiling then we will suppose that a family
of tiles is (arbitrarily) distinguished, and we will call it the
$D$-th family. We will also denote
by $\bar{t}$ the tiling obtained from $t$ when we delete the $D$-th
family of tiles. Therefore, $\bar{t}$ is a \ran{D-1}{d}\ tiling.

\fig{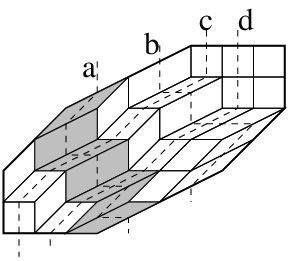}{Some de Bruijn lines of a two dimensional tiling.
Lines $a$ and $b$ belong to the same family, whereas they are
not in the same family as $c$ and $d$. Moreover, $c$ and $d$
belong to the same family.
The shaded tiles represent the worm associated to the line $a$.}

In the following, we will use a variant of the classical
de Bruijin lines: the so-called \emph{oriented} de Bruijn lines,
which are simply the de Bruijn lines together with an orientation
over each line such that for all family, each line in this family
have the same orientation. This lead to the usual notion of
\emph{dual graph} of a tiling: its set of vertices is the
set of intersection points of de Bruijn lines, and there is an
edge $(i,j)$ if and only if $i$ and $j$ are in adjacent tiles
and there is a piece of de Bruijn line oriented from $a$ to $b$.
Notice that dual graphs of tilings are usually undirected graphs
representing the neighbourhood relation of the tiles in the tiling.
The two definitions are equivalent, except the orientation of
the edges, which is necessary in the following.

Despite our results may be general, we will restrict ourselves
to the $d=2$ special case in the following. This means that we tile
two dimensional zonotopes (i.e. $2D$-gons) with lozenges.
This restriction is due
to the fact that these tilings received most of the attention until
now, which allows us to use some previously known results which
have not yet been established in the general $d$ dimensional case
(and which may be false in this case). Notice however that
an isomorphism between \ran{d+1}{d}\ tilings, ordered with
the transitive and reflexive closure of the flip relation,
and $d$-dimensional hypersolid partitions ordered with
the transitive and reflexive closure of the addition of one
grain, is exhibited in \cite{DMB97}. Therefore, we can already
say from Theorem~\ref{th_part} that the set of \ran{d+1}{d}\ tilings
is a distributive lattice.

In \cite{Des97,DMB97,DMB99}, Destainville studied the relation between
rhombic tilings of zonotopes and integer partitions.
In the following, we will widely use the correspondence he exhibited.
Therefore, we shortly describe it here.
For more details, we refer to the original papers.

First, let us see
how he associates a partition to a tiling. Let $t$ be a tiling of a
zonotope $Z$. Let us consider $\bar{t}$,
the tiling obtained from $t$ by deleting the tiles of the $D$-th family.
Now consider the oriented de Bruijn lines of $t$.
For any tile $\tau$ which is not in the $D$-th family (i.e. a tile
in $\bar{t}$), and any de Bruijn line $l$ which is not in the $D$-th
family, we define
$w_{\tau,l}$ as the number of de Bruijn lines of the $D$-th family
we cross if we start
from the tile $\tau$ and follow the de Bruijn line $l$
(with respect to its orientation). Notice that, since the de Bruijn
lines of the $D$-th family define a partition of the tiling
into disjoint regions which can not touch two opposite borders of $Z$,
and since the de Bruijn lines go from one border of $Z$ to its opposite,
we can always choose the orientations to have $w_{\tau,l_1} = w_{\tau,l_2}$
where $l_1$ and $l_2$ are the two de Bruijn lines which cross $\tau$.
Therefore, we can denote
this value by $w_{\tau}$. Now, consider the dual graph of $\bar{t}$,
$G=(V,E)$. Then, the function $p$ defined for all $v$ in $V$ by $p(v) =
w_{\tau}$, where $\tau$ is the dual tile of the vertex $v$, is a
partition solution to the partition problem $(G,h)$ where $h$ is the
total number of de Bruijn lines in the $D$-th family in $t$. In the
following, given a tiling $t$, we will denote by $\mathcal{P}(t)$ the partition
associated this way to $t$.

\begin{figure}
\centerline{\input{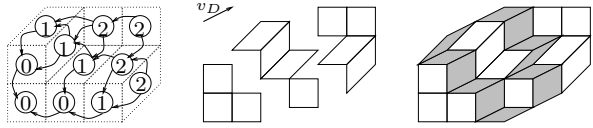}}
\caption{\label{fig_T}From left to right : A partition which is solution
of a partition problem over a graph which is the dual graph of a
\ran{3}{2}\ tiling (dotted), the tiles of this tiling translated with
respect to the values of the corresponding parts and the vector $v_D$,
and the obtained \ran{4}{2}\ tiling after completion. This gives an
example of the function $\mathcal{T}$. Notice that we can do the way
back from right to left and obtain this way an example for $\mathcal{P}$.}
\end{figure}

Conversely, given a partition $p$ solution of the partition
problem $(G,h)$ where $G$ is the dual graph of a tiling $t$, we want
to define a tiling $t'$ associated to $p$.
Let $Z$ be the zonotope tiled by $t$. Let $Z'$ be the zonotope
generated by the same
family of vectors than $Z$ with an additional one: $v_D$ with $l_D = h$.
Let us consider the following partition of the set of vertices of $G$
(and dually of the tiles of $t$) : $V_i = \lbrace v \in V \mbox{ such that }
p_v = i \rbrace$. We can now construct $t'$ by insertion of a worm of
the $D$-th family in $t$ between the tiles corresponding to the
sets $V_i$ and $V_{i+1}$ for all $i$. In other words, starting from $t$,
for all $i$ we apply the translation defined by the vector
$i\cdot v_D$ to the
tiles in $V_i$, and we add the tiles of the $D$-th family in order
to fill $Z'$. We obtain this way $t'$, the tiling of $Z'$ associated
to $p$. In the following, given a partition $p$, we will denote by
$\mathcal{T}(p)$ the partition associated this way to $p$. See
Figure~\ref{fig_T} for an example.

Destainville obtains this way a method to generate all the \ran{D}{2}\ 
tilings of a given zonotope $Z$. He starts from a \ran{3}{2}\ tiling,
which is nothing but the projection of the Ferrer's diagram of a
planar partition \cite{DMB97}. Indeed, a planar partition can be
viewed as a decreasing sequence of stacked cubes on a $2$-dimensional
grid, the projection of which gives the \ran{3}{2}\ tiling. See
Figure~\ref{fig_ex_til} (left) for an example. From the oriented dual graph
of this tiling, he defines a partition problem, the solutions of which
are equivalent to \ran{4}{2}\ tilings, as explained above. Likewise,
he can construct a \ran{D+1}{2}\ tiling from a \ran{D}{2}\ one for
any $D$, and so obtains a way to generate \ran{D}{2}\ tilings for
any $D$. It is shown in \cite{Des97,DMB99} that the application
$\mathcal{T}$, which generate the tiling associated to
a partition is a bijection from the set of the partitions solutions
to the problems $(G_t,h)$ for all oriented dual graph $G_t$ of a
\ran{D}{d}\ tiling to the set of \ran{D+1}{d}\ tilings. Moreover,
it is shown in these papers that this bijection is an order
isomorphism if the set of partitions is ordered with the reflexive
and transitive closure of the relation induced by the additions of
one grain (see Section~\ref{sec_part}), and if the set of tilings
is ordered with the reflexive and transitive closure of the
relation induced by the flips.

It was shown
in \cite{Ken93} and \cite{Eln97} that 
we can obtain all the two dimensional tilings of
a zonotope from a given one by iterating the following rule:
the transition $t \rightarrow t'$ is possible if the tiling $t'$ can
be obtained from $t$ with a flip\,\footnote{The generalisation of
this claim to any dimension, however, has been pointed out by
Reiner as a difficult open question \cite{Rei99}.}.
This transition rule leads us to consider the tilings of a zonotope as
the possible states of a dynamical system.
A sequence of
such transitions is denoted by $\flips$, which is equivalent
to the transitive and reflexive closure of $\rightarrow$, also denoted
by $\ge$, depending on the emphasis given to the dynamical aspect
or to the order theoretical approach.
We denote by $T(Z,D,2)$ the set of \ran{D}{2}\ tilings
of the zonotope $Z$ ordered by $\ge$.
An example is given in Figure~\ref{fig_pav_tr}.
Notice that all sequences of flips (with no inverse flips) from
a tiling to another one have the same length \cite{DMB99}, which will
be useful in the following.
Using the preliminaries given above, we can now state the first results we
need to prove that $T(Z,D,2)$ is a lattice.

\begin{lemme}
\label{lem_struct}
The set $T(Z,D,2)$ is the disjoint union of distributive
lattices $L_i$ such that
a flip transforms a tiling in $L_i$ into another one in $L_i$ if and only
if it involves at least one tile of the $D$-th family.
Moreover, for all $t \in L_i$ and $u \in L_j$,
$\bar{t} = \bar{u} \Leftrightarrow i = j$.
\end{lemme}
\vspace*{-0.1cm}
\proof
Let us consider the maximal subsets $L_i$ of $T(Z,D,2)$ such that
a flip goes from a tiling in $L_i$ to another one in $L_i$ if and only
if it involves at least one tile of the $D$-th family.
It is shown in \cite{Des97} and \cite{DMB99} that such a set, equipped
with the transition rule described above (flip), is isomorphic to the
set of the solutions of a partition problem, depending on
$Z$ and $D$, equipped with the transition rule described in
Section~\ref{sec_part} (addition of one grain). We know from
Theorem~\ref{th_part} that this set is a distributive lattice.
Therefore, we obtain the first part of the claim.\\
It is then clear that if $s$ and $t$ are in $L_i$ then $\bar{s}=\bar{t}$:
it suffices to notice that if $t \longrightarrow t'$ such that this flip
involves at least one tile in the $D$-th family then $\bar{t}=\bar{t'}$.
Moreover if $\bar{s}=\bar{t}$ then $s$ can not be obtained from $t$ with
a flip involving three tiles with none of them belonging to the $D$-th family: such a flip
changes the position of the tiles in $\bar{s}$ and $\bar{t}$.
This ends the proof.
\cqfd

\begin{lemme}
\label{lem_max_min}
Let $a$, $b$ and $c$ be in a $L_i$ ($L_i$ being one of the sets partitioning
$T(Z,D,2)$ defined in Lemma~\ref{lem_struct})
such that $a$ is the unique maximal element of $L_i$ and $b$
is its unique minimal element.
If a flip involving three tiles none of them belonging to
the $D$-th family is possible from $c$ then it is possible from $a$ and $b$.
\end{lemme}
\vspace*{-0.1cm}
\proof
First notice that, since a flip inside $L_i$ involves tiles which
are in the $D$-th family, we have $\bar{a}=\bar{b}=\bar{c}$.
Therefore, the flip from $c$ is possible from $a$ and $b$
if the three tiles it involves are neighbours in $a$ and $b$.
From Lemma~\ref{lem_struct}, $\mathcal{P}(a)$ and $\mathcal{P}(b)$,
the partitions which correspond to the tilings $a$ and $b$, are 
respectively the maximal and minimal elements of the set $P(G,h)$
of solutions to a partition problem $(G,h)$. Therefore, from
Theorem~\ref{th_part}, $\mathcal{P}(a)$ and $\mathcal{P}(b)$
have all their parts equal to respectively $0$ and $h$.
Then, from the definition of $\mathcal{T}=\mathcal{P}^{-1}$,
all the tiles which do not belong to the $D$-th family tile
a sub-zonotope of $Z$, and so all the flips involving
three such tiles are possible from $a$ and $b$.
\cqfd

\fig{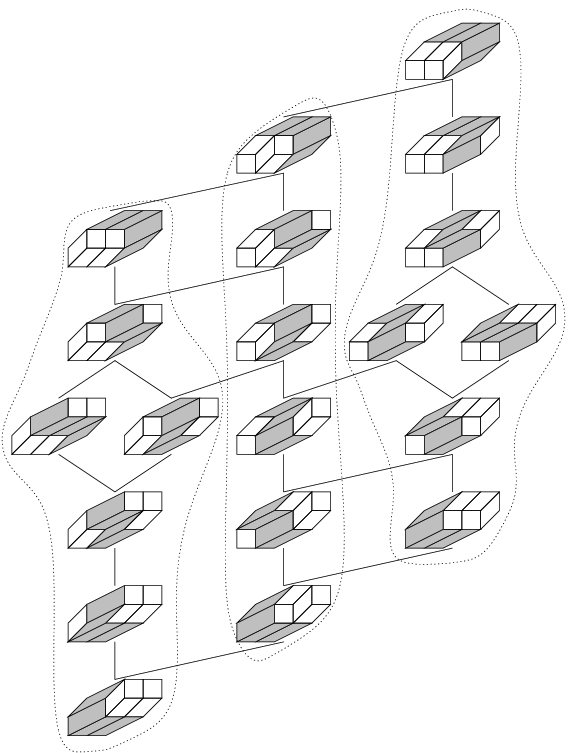}{$T(Z,4,2)$ for a given $Z$. The possible transitions (flips)
are represented. The shaded tiles show the $4$-th family, and the
dotted sets are the distributive lattices $L_1$, $L_2$ and $L_3$,
as stated by Lemma~\ref{lem_struct}.}

We can now define the order $\overline{T(Z,D,2)}$ as the quotient of
$T(Z,D,2)$ with respect to the equivalence relation
$s \sim t \Leftrightarrow \bar{s} = \bar{t}$, i.e. defined by
the lattices $L_i$. In other words,
we consider the set of the lattices $L_i$ as the set of vertices
of $\overline{T(Z,D,2)}$, and there is one edge from $L_i$ to $L_j$
in $\overline{T(Z,D,2)}$ if and only if there is at least a transition
from one element of $L_i$ to one element of $L_j$ in $T(Z,D,2)$.
Also notice that if, given a tiling of a zonotope Z, we delete the tiles
in the $D$-th family, we obtain a new zonotope. Notice that this zonotope
only depends on $Z$ and does not depend of the considered tiling, since,
as one can easily verify, if $t \longrightarrow t'$ then $\bar{t}$ and
$\bar{t'}$ tile the same zonotope.
Let us denote by $Z'$ this zonotope.

\begin{lemme}
\label{lem_quot}
The order $\overline{T(Z,D,2)}$ is isomorphic to the order $T(Z',D-1,2)$.
\end{lemme}
\vspace*{-0.1cm}
\proof
From Lemma~\ref{lem_struct}, we can associate to
each $L_i \in \overline{T(Z,D,2)}$
a tiling $t_i$ such that for all $a$ in $L_i$, $\bar{a}=t_i$.
It is clear that $t_i$ is in $T(Z',D-1,2)$.
Conversely, if we have a tiling of $Z'$, then we can use the construction
of a \ran{D+1}{d}\ tiling from a \ran{D}{d}\ one described above to
obtain a tiling $t$ of $Z$. Therefore, there is a bijection
between $\overline{T(Z,D,2)}$ and $T(Z',D-1,2)$.
We will now see that it is an order isomorphism.
From Lemma~\ref{lem_struct}, if there exists a flip
$a \rightarrow b$ between two tilings $a\in L_i$ and $b\in L_j$ with
$i \not= j$, then it does not involve any tile of the $D$-th family,
and so there exists a flip $t_i \rightarrow t_j$ in
$T(Z',D-1,2)$. Conversly, if there is a flip $t_i \rightarrow t_j$
in $T(Z',D-1,2)$, then there exists $a \in L_i$ and $b \in L_j$ such
that $a \rightarrow b$: from Lemma~\ref{lem_max_min}, it
suffices for example to take the maximal
elements of $L_i$ and $L_j$ respectively for $a$ and $b$.
Therefore, $\overline{T(Z,D,2)}$ is isomorphic to $T(Z',D-1,2)$.
\cqfd

With these two lemmas, we have much information about any set
$T(Z,D,2)$: it is the disjoint union of distributive lattices,
and its quotient with respect to these lattices has itself the
same structure, since it is isomorphic to $T(Z',D-1,2)$. This shows
that the sets $T(Z,D,2)$ are strongly structured, and makes
it possible to write efficient algorithms based on this structure,
for example the computation of a shortest sequence of flips which
transforms a given tiling into another given one.

\section{Conclusion and perspectives}

In conclusion, we gave structural results on generalized integer
partitions and two dimensional tilings of zonotopes. Our main tools
were dynamical systems
and order theories. This allows an intuitive presentation of
the topic, and the presence of lattices in this kind of dynamical
systems seems very general \cite{Pha99}. This makes it possible to develop
efficient algorithms for a variety of questions over tilings. We wrote
for example an algorithm which transforms a tiling into another
one with a minimal number of flips (we do not give it here because
of the lack of space).
The lattice structure is also strongly related to enumeration
problems, and often
offers the possibility of giving new combinatorial results. For
example, one could use the results presented here to look for
a formula for the minimum number of flips necessary to transform a
tiling into another one.

There are two immediate directions in which it seems promising to extend the
results presented here. The first and obvious one is to study $d$ dimensional
tilings with $d>2$. Then, it is not clear wether all the tilings of a given
zonotope can be obtained from a particular one by flipping tiles \cite{Rei99}.
Moreover, we would need an efficient formalism in order to give clear and
nice proofs. This formalism still has to be developed. The other important
remark is that the choice of the $D$-th family all along our work is
arbitrary. This means that we could choose any family to be the $D$-th,
and so there are many ways to decompose $T(Z,D,2)$ into a disjoint union of
distributive lattices. This is a strong and surprinsing fact, which has
to be fully explored.

Notice also that Remila \cite{Rem99a, Rem99b} showed that special classes of
tilings with flips are lattices: the domino tilings, the bar tilings
and the calisson tilings. Notice that this last class is nothing but
the \ran{3}{2}\ tilings. However, the techniques he used to
prove these results are very different from the one presented here, and
it would be very interesting to try and extend these results to
larger classes of tilings. The idea of considering tilings as projections
of some high-dimensional structures seems promising, since it allowed
us to study the large class of tilings. Moreover, the bars
and dominoes tilings studied by Remila can also be viewed as
projections of high-dimensional objects. Therefore, this approach
would be an interesting research direction for general results.

Finally, one may wonder if the results presented here always stands
when the support of the tiling is not a zonotope. It would be
interesting to know the limits of our structural results. They may
be very general, and lead to sub-lattices properties of the obtained
sets of tilings. Likewise,
it would be useful to study what happens when the size of the
zonotope grows to infinity. Some results about that are presented in
\cite{Des97} and \cite{DMB99} but a lot of work remains to be done.

\section{Erratum}
In a previous version of this paper, we claimed that the sets
$T(Z,D,2)$ were themselves lattices, which is actually false.
Indeed, the set $T(Z,D,2)$ when $Z$ is a $2d$-gon having $l_i=1$
for all $i$ is isomorphic to the higher bruhat order $B(n,2)$
(see \cite{FZ99}), and it is known that $B(6,2)$ is \emph{not}
a lattice (\cite{Zie93}, Theorem 4.4). We apologize for this,
and we thank V.~Reiner, who first pointed out this error and
gave us useful references.

\small
\section{Acknowledgments}
I thank Nicolas Destainville, Michel Morvan, Eric Remila and Laurent
Vuillon for many useful comments on preliminary versions. Special
thanks are due to Laurent Vuillon who introduced me to the
facinating field of tilings. I also thank Serge Elnitsky, Richard
Kenyon, Vic Reiner, G\"unter Ziegler and the dominoers for helpful
references.

\bibliographystyle{plain}
\bibliography{../Bib/bib}

\end{document}